\pdfoutput=1
\RequirePackage{ifpdf}
\ifpdf % We are running pdfTeX in pdf mode
\documentclass[pdftex]{sigma}
\else
\documentclass{sigma}
\fi

\numberwithin{equation}{section}

\begin{document}

\newcommand{\arXivNumber}{1310.2006}

\allowdisplaybreaks

\renewcommand{\PaperNumber}{069}

\FirstPageHeading

\ShortArticleName{Special Solutions and Linear Monodromy for G(1112)}

\ArticleName{Special Solutions and Linear Monodromy\\
for the Two-Dimensional Degenerate Garnier\\
System G(1112)}

\Author{Kazuo KANEKO}

\AuthorNameForHeading{K.~Kaneko}

\Address{Seki Kowa Institute of Mathematics, Yokkaichi University,\\
Kayaucho, Yokkaichi, Mie, 512-8512, Japan}
\Email{\href{mailto:dr_kaneko_k@yahoo.co.jp}{dr\_kaneko\_k@yahoo.co.jp}}

\ArticleDates{Received October 24, 2013, in f\/inal form June 14, 2014; Published online July 05, 2014}

\Abstract{We have classif\/ied special solutions around the origin for the two-dimensional degenerate Garnier system
G(1112) with generic values of complex parameters, whose linear monodromy can be calculated explicitly.}


\Classification{34M55; 33C15}

\section{Introduction}
We have studied special solutions with generic values of complex parameters for the fourth, f\/ifth, sixth and third
Painlev\'e equations, for which the monodromy data of the associated linear equation (we call linear monodromy) can be
calculated explicitly~\cite{KK1,KK2,KK4,KO}.
These papers are based on A.V.~Kitaev's idea who calculated f\/irst the linear monodromy with generic value of complex
parameter explicitly by taking examples of the f\/irst and second Painlev\'e equations~\cite{AVK}.
We remark that P.~Appell~\cite{AP} also studied the symmetric solutions to the f\/irst and second Painlev\'e equations,
but he did not study linear monodromy problems.

The Garnier system was derived by R.~Garnier (1912) as the extension of the sixth Painlev\'e equation~\cite{Gar}.
The original Garnier system has~$n$ variables and is
expressed in the nonlinear partial dif\/ferential equations system,
whose dimension of the solution space is $2n$.
There are few research for the special solutions to the Garnier system compared with Painlev\'e equations.
We will study the Garnier transcendents by applying f\/irst
the same method to the two-dimensional Garnier system,
which we have used for the Painlev\'e equations above.
Some new discovery is expected by
viewing
Painlev\'e equations from the Garnier system.

Two-dimensional Garnier system has the following degeneration diagram similar to the Pain\-le\-v\'e
equations~\cite{Kim}:
\begin{gather*}
\begin{matrix}
{\rm G}(11111)& \rightarrow & {\rm G}(1112) & \rightarrow & {\rm G}(122) & \rightarrow & {\rm G}(23) & &
\\
& & \downarrow & & \downarrow & & \downarrow & &
\\
& & {\rm G}(113)& \rightarrow & {\rm G}(14) & \rightarrow & {\rm G}(5)& \rightarrow & {\rm G}(9/2)&
\end{matrix}
\end{gather*}
(The degeneration from {\rm G}(113) to {\rm G}(23) also exists.) Numbers in
brackets
represents a~partition of~5. The number~1 represents the regular singular point
 and the number $r+1$ represents an irregular singular point of Poincar\'e rank~$r$.
The two-dimensional Garnier system {\rm G}(11111) which is the extension of the sixth Painlev\'e equation ${\rm P}_{\rm VI}$
degenerates step by step to the two-dimensional degenerate Garnier system {\rm G}(9/2) which is the extension of the f\/irst
Painlev\'e equation ${\rm P}_{\rm I}$.

The purpose of this paper is to obtain the special solutions to the system {\rm G}(1112), for which the linear monodromy
$\big\{M_0,M_1=S_1^{(1)}S_2^{(1)}e^{2\pi i T_1},M_{t_2},M_{\infty}\big\}$ can be calculated explicitly.

The two-dimensional degenerate Garnier system {\rm G}(1112)~$\{K_1, K_2, \lambda_1, \lambda_2, \mu_1, \mu_2, t_1, t_2\}$
is derived as the extension of the f\/ifth Painlev\'e equation by the isomonodromic deformation of the second kind, non-Fuchsian
ordinary dif\/ferential equation, which has three regular singularities and one irregular singularity of
Poincar\'e rank~1 on the Riemann sphere~\cite{Kim, Ok}:
\begin{gather}
\frac{d^2\psi}{dx^2}+\biggl[\frac{1-\alpha_0}{x}+\frac{\eta t_1}{(x-1)^2}+\frac{2-\alpha_1}{x-1}
+\frac{1-\alpha_2}{x-t_2} -\frac{1}{x-\lambda_1}-\frac{1}{x-\lambda_2}\biggr]\frac{d\psi}{dx}
\nonumber
\\
\hphantom{\frac{d^2\psi}{dx^2}}
{}+\biggl[\frac{\nu(\nu+\alpha_{\infty})}{x(x-1)}-\frac{t_1K_1}{x(x-1)^2} -\frac{t_2(t_2-1)K_2}{x(x-1)(x-t_2)}
\nonumber
\\
\hphantom{\frac{d^2\psi}{dx^2}}
{}+\frac{\lambda_1(\lambda_1-1)\mu_1}{x(x-1)(x-\lambda_1)}
+\frac{\lambda_2(\lambda_2-1)\mu_2}{x(x-1)(x-\lambda_2)}\biggr]\psi=0,
\label{lg2}
\end{gather}
where $K_1$ and $K_2$ are Hamiltonians, $\lambda_1$, $\lambda_2$, $\mu_1$ and $\mu_2$ are the Garnier functions, $t_1$
and $t_2$ are deformation parameters and $\alpha_j$ $(j=0,1,2,\infty)$, $\nu \in \mathbb{C}$ and $\eta \in \mathbb{C^{\times}}$ are
complex parameters.
The Riemann scheme of~\eqref{lg2} is
\begin{gather*}
P\left(
\begin{matrix}
x=0 & x=1 &x=t_2 &x=\lambda_1 &x=\lambda_2 & x= \infty &
\\
\begin{matrix}
0
\\
\alpha_0
\end{matrix}
&
\overbrace{
\begin{matrix}
0 & 0
\\
\eta t_1 & \alpha_1
\end{matrix}
}
&
\begin{matrix}
0
\\
\alpha_2
\end{matrix}
&
\begin{matrix}
0
\\
2
\end{matrix}
&
\begin{matrix}
0
\\
2
\end{matrix}
&
\begin{matrix}
\nu
\\
\nu+\alpha_{\infty}
\end{matrix}
\begin{matrix}
;x
\end{matrix}
\end{matrix}
\right),
\\
\alpha_0+\alpha_1+\alpha_2+\alpha_{\infty}=1-2\nu.
\end{gather*}

This is also derived by the conf\/luence of two regular singularities $x=t_1$ and $x=1$ in the two-dimensional Garnier system {\rm G}(11111).

{\rm G}(1112) has movable algebraic branch points and Hamiltonian structure expressed in rational function.
We have the two-dimensional degenerate Garnier system $\mathcal{H}_2(1112)\{H_1,H_2,q_1,q_2$, $p_1,p_2,s_1,s_2\}$  by the canonical transformations:
\begin{gather*}
 s_1=\frac{1}{t_1},
\qquad
s_2=\frac{t_2}{t_2-1},
\qquad
-t_1(t_2-1)q_1=(\lambda_1-1)(\lambda_2-1),
\\
(t_2-1)^2q_2=(\lambda_1-t_2)(\lambda_2-t_2),
\qquad
\mu_i=\frac{q_1p_1}{\lambda_i-1}+\frac{q_2p_2}{\lambda_i-t_2},
\qquad
i=1, 2,
\\
 H_1=-t_1^2\left(K_1+\sum\limits_{j=1}^{2}p_j\frac{\partial q_j}{\partial t_1}\right),
\qquad
H_2=-(t_2-1)^2\left(K_2+\sum\limits_{j=1}^{2}p_j\frac{\partial q_j}{\partial t_2}\right),
\\
 \sum\limits_{j=1}^{2}(dp_j\wedge dq_j-dH_j\wedge ds_j)=\sum\limits_{j=1}^{2}(d\mu_j\wedge d\lambda_j-dK_j\wedge dt_j).
\end{gather*}

$\mathcal{H}_2$ has the Painlev\'e property and the polynomial Hamiltonian structure~\cite{IKSY,Kim, Ok}.
We obtain the special solutions in the Hamiltonian system $\mathcal{H}_2$
and then inversely transform them to the solutions in the Hamiltonian system {\rm G}(1112), which  are
substituted into the linear equation~\eqref{lg2}.
We obtain
eight meromorphic solutions with generic values of complex parameters around the origin $(t_1, t_2)=(0,0)$, which we
name the solutions (1), (2), \dots, (8).

The calculation of the linear monodromy consists of three steps.
The f\/irst step is taking the limit $(t_1, t_2)\rightarrow (0,0)$ after substituting the solution into the linear
equation~\eqref{lg2}.
We call this step ``the f\/irst limit'', in which the linear monodromy matrices
$M_0$ and $M_{t_2}$ are calculated as the conf\/luent
linear monodromy~$M_{t_2}M_0$.
In the second step, we separate this conf\/luent linear monodromy~$M_{t_2}M_0$.
After transforming the linear equation~\eqref{lg2} by putting $x=t_2\xi$ and substituting the solution into the linear
equation~\eqref{lg2}, we take the limit $(t_1, t_2)\rightarrow (0,0)$.
We call this step ``the second limit''.
In the third step, we transform the linear equation~\eqref{lg2} by putting $x-1=\eta t_1/z$ which keeps the
irregularity at $x=1$ so that we can calculate the Stokes matrices $\big\{S_1^{(1)},S_2^{(1)}\big\}$.
We call this step ``the third limit''.

Each of the obtained eight meromorphic solutions with generic values of complex parameters around the origin $(t_1, t_2)=(0,0)$
has the remarkable characteristics, respectively.
The four solutions make the two monodromy matrices commutable and the Stokes matrices around $x=1$ unity and the other
four solutions make the three monodromy matrices commutable, which are summarized in Theorem~\ref{tm:4}.

In Appendix~\ref{appendix1}, we show the fundamental solutions and the associated monodromy matrices of Gauss
hypergeometric equation and Kummer's equation.
In Appendix~\ref{appendix2}, we show the Briot--Bouquet's theorem for a~system
of partial dif\/ferential equations in two variables and short comment on it, how it proves convergence of the eight solutions.

\section[The two-dimensional degenerate Garnier system $\mathcal{H}_2(1112)$]{The two-dimensional degenerate
Garnier system $\boldsymbol{\mathcal{H}_2(1112)}$}

In this section, we write down the polynomial Hamiltonians $H_1$, $H_2$ and the Hamiltonian system $\mathcal{H}_2(1112)$.

$\bullet$ Hamiltonians $H_1$ and $H_2$:
\begin{gather*}
s_1^2H_1 = q_1^2(q_1-s_1)p_1^2+2q_1^2q_2p_1p_2+q_1q_2(q_2-s_2)p_2^2
\\
\phantom{s_1^2H_1 =}
{}- \biggl[(\alpha_0+\alpha_2-1)q_1^2+\alpha_1q_1(q_1-s_1)+\eta(q_1-s_1)+\eta s_1q_2\biggr]p_1
\\
\phantom{s_1^2H_1 =}
{}- \biggl[(\alpha_0+\alpha_1-1)q_1q_2+\alpha_2q_1(q_2-s_2)-\eta(s_2-1)q_2\biggr]p_2
+ \nu(\nu+\alpha_{\infty})q_1,
\\
s_2(s_2-1)H_2 = q_1^2q_2p_1^2+2q_1q_2(q_2-s_2)p_1p_2
\\
\phantom{s_1^2H_1 =}
{}+ \biggl[q_2(q_2-1)(q_2-s_2)+\frac{s_2(s_2-1)}{s_1}q_1q_2\biggr]p_2^2
\\
\phantom{s_1^2H_1 =}
{}- \biggl[(\alpha_0+\alpha_1-1)q_1q_2+\alpha_2q_1(q_2-s_2)-\eta(s_2-1)q_2\biggr]p_1
\\
\phantom{s_1^2H_1 =}
{}- \biggl[(\alpha_0-1)q_2(q_2-1)+\alpha_1q_2(q_2-s_2)+\alpha_2(q_2-1)(q_2-s_2)
\\
\phantom{s_1^2H_1 =}
{}+\frac{s_2(s_2-1)}{s_1}(\alpha_2q_1+\eta q_2)\biggr]p_2+\nu(\nu+\alpha_{\infty})q_2.
\end{gather*}

$\bullet$ Hamiltonian system $\mathcal{H}_2(1112)$:
\begin{gather*}
-t_1\frac{\partial q_1}{\partial t_1}=s_1\frac{\partial H_1}{\partial p_1},
\qquad
-t_1\frac{\partial q_2}{\partial t_1}=s_1\frac{\partial H_1}{\partial p_2},
\qquad
t_1\frac{\partial p_1}{\partial t_1}=s_1\frac{\partial H_1}{\partial q_1},
\qquad
t_1\frac{\partial p_2}{\partial t_1}=s_1\frac{\partial H_1}{\partial p_1},
\\
\frac{\partial q_1}{\partial s_2}=\frac{\partial H_2}{\partial p_1},
\qquad
\frac{\partial q_2}{\partial s_2}=\frac{\partial H_2}{\partial p_2},
\qquad
-\frac{\partial p_1}{\partial s_2}=\frac{\partial H_2}{\partial q_1},
\qquad
-\frac{\partial p_2}{\partial s_2}=\frac{\partial H_2}{\partial q_2}.
\end{gather*}

\begin{remark}
We use $t_1$ $(=1/s_1)$
instead of $s_1$ to apply the Briot--Bouquet's theorem~\cite{BB} at the origin.
\end{remark}

\section[Meromorphic solutions around the origin $(t_1=1/s_1,s_2)=(0,0)$]{Meromorphic solutions\\
around the origin $\boldsymbol{(t_1=1/s_1,s_2)=(0,0)}$}

In this section, we give the calculated meromorphic solutions around $(t_1=1/s_1, s_2) =(0, 0)$, which are satisf\/ied
with the Hamiltonian system $\mathcal{H}_2$(1112).

When $q_i$ and $p_i$ $(i=1,2)$ are meromorphic, they have at most a~simple pole around $(t_1=1/s_1, s_2) =(0, 0)$.
Let $R=\mathbb{C}\{\{t_1, s_2\}\}$ be the ring of convergent power series of $t_1$ $(=1/s_1)$ and~$s_2$ around the origin
and for $u_1, u_2, \dots,u_n \in R$, let $\langle u_1, u_2, \dots,u_n \rangle$ be the ideal of~$R$ generated by $u_1,
u_2,\dots, u_n$.
\begin{theorem}
For generic values of the complex parameters $\{\alpha_0, \alpha_1, \alpha_2, \alpha_{\infty},\nu,\eta \}$, the
Hamiltonian system $\mathcal{H}_2(1112)$ has the following eight meromorphic solutions around $(t_1=1/s_1, s_2) =(0,0)$:
\begin{gather*}
(1)\ \
q_1=\frac{\eta}{\alpha_{\infty}}+\langle t_1,s_2\rangle,
\qquad
q_2=\frac{\alpha_{\infty}+\alpha_1}{\alpha_{\infty}} +\langle t_1,s_2\rangle,
\\
\phantom{(1)}\ \
p_1=\frac{(\alpha_{\infty}+\alpha_1)(\nu+\alpha_2)}{\alpha_{\infty}(1-\alpha_{\infty}-\alpha_1)}t_1s_2
+\langle t_1^3,t_1^2s_2, s_2^2\rangle,
\qquad
p_2=\frac{-\nu\alpha_{\infty}}{\alpha_{\infty}+\alpha_1}+\langle t_1,s_2\rangle,
\\
(2)\ \
q_1=\frac{-\eta}{\alpha_{\infty}}+\langle t_1,s_2\rangle,
\qquad
q_2=\frac{\alpha_{\infty}-\alpha_1}{\alpha_{\infty}} +\langle t_1,s_2\rangle,
\\
\phantom{(2)}\ \
p_1=\frac{(\alpha_{\infty}-\alpha_1)(\nu+\alpha_2+\alpha_{\infty})}
{\alpha_{\infty}(1+\alpha_{\infty}-\alpha_1)}t_1s_2
+\langle t_1^3,t_1^2s_2, s_2^2\rangle,
\qquad
p_2=\frac{-(\alpha_{\infty}+\nu)}{\alpha_{\infty}-\alpha_1}\alpha_{\infty}
+\langle t_1,s_2\rangle,
\\
(3)\ \
q_1=\frac{-\eta}{\alpha_1}+\langle t_1,s_2\rangle,
\qquad
q_2=s_2\left(\frac{\alpha_2}{\alpha_0+\alpha_2} +\langle t_1,s_2\rangle \right),
\\
\phantom{(3)}\ \
p_1=\frac{\nu(\nu+\alpha_{\infty})}{1-\alpha_1}t_1+\langle s_2,t_1^2\rangle,
\qquad
p_2=\frac{1}{s_2}\left(\frac{\nu(\nu+\alpha_{\infty})}{1-\alpha_0-\alpha_2}s_2+\langle t_1^2,t_1s_2, s_2^2\rangle
\right),
\\
(4)\ \
q_1=\frac{-\eta}{\alpha_1}+\langle t_1,s_2\rangle,
\qquad
q_2=s_2\left(\frac{\alpha_2}{\alpha_2-\alpha_0} +\langle t_1,s_2\rangle\right),
\\
\phantom{(4)}\ \
p_1=\frac{(\nu+\alpha_2)(\nu+\alpha_2+\alpha_{\infty})}{1-\alpha_1}t_1+\langle s_2,t_1^2\rangle,
\qquad
p_2=\frac{1}{s_2}\left(\alpha_0-\alpha_2+\langle t_1,s_2\rangle\right),
\\
(5)\ \
q_1=\frac{1}{t_1}\left(\frac{\alpha_{\infty}+\alpha_0+\alpha_2}{\alpha_{\infty}}+\langle t_1,s_2\rangle\right),
\qquad
q_2=\frac{-\alpha_2}{\alpha_{\infty}}s_2+\langle t_1^2,t_1s_2,s_2^2\rangle,
\\
\phantom{(5)}\ \
p_1=t_1\left(\frac{-\nu\alpha_{\infty}}{\alpha_{\infty}+\alpha_0+\alpha_2}+\langle t_1,s_2\rangle\right),
\qquad
p_2=\frac{\eta\nu\alpha_{\infty}t_1}{(\alpha_1+2\nu)(\alpha_{\infty}+\alpha_0+\alpha_2)}+\langle s_2,t_1^2\rangle,
\\
(6)\ \
q_1=\frac{1}{t_1}\left(\frac{\alpha_{\infty}-\alpha_0-\alpha_2}{\alpha_{\infty}}+\langle t_1,s_2\rangle\right),
\qquad
q_2=\frac{\alpha_2}{\alpha_{\infty}}s_2+\langle t_1^2,t_1s_2, s_2^2\rangle,
\\
\phantom{(6)}\ \
p_1 = t_1\left(\frac{-(\nu+\alpha_{\infty})\alpha_{\infty}}{\alpha_{\infty}-\alpha_0-\alpha_2}+\langle t_1,s_2\rangle
\right),
\\
\phantom{(6)}\ \
p_2=\frac{\eta\alpha_{\infty}(\nu+\alpha_{\infty})t_1}{(\alpha_{\infty}-\alpha_0-\alpha_2)(1-\alpha_{\infty}+\alpha_0+\alpha_2)}+\langle
s_2,t_1^2\rangle,
\\
(7)\ \
q_1=\frac{1}{t_1}\left(\frac{\alpha_{\infty}+\alpha_0-\alpha_2}{\alpha_{\infty}}+\langle t_1,s_2\rangle\right),
\qquad
q_2=s_2\left(\frac{\alpha_2}{\alpha_{\infty}}+\langle t_1,s_2\rangle\right),
\\
\phantom{(7)}\ \
p_1=t_1\left(\frac{-(\nu+\alpha_2)\alpha_{\infty}}{\alpha_{\infty}+\alpha_0-\alpha_2}+\langle t_1,s_2\rangle\right),
\qquad
p_2=\frac{1}{s_2}\left(\alpha_{\infty}+\langle t_1,s_2\rangle\right),
\\
(8)\ \
q_1=\frac{1}{t_1}\left(\frac{\alpha_{\infty}-\alpha_0+\alpha_2}{\alpha_{\infty}}+\langle t_1,s_2\rangle\right),
\qquad
q_2=s_2\left(\frac{-\alpha_2}{\alpha_{\infty}}+\langle t_1,s_2\rangle\right),
\\
\phantom{(8)}\ \
p_1=t_1\left(\frac{-(\nu+\alpha_2+\alpha_{\infty})\alpha_{\infty}}{\alpha_{\infty}-\alpha_0+\alpha_2} +\langle
t_1,s_2\rangle\right),
\qquad
p_2=\frac{1}{s_2}\left(-\alpha_{\infty}+\langle t_1,s_2\rangle\right).
\end{gather*}
\end{theorem}
\begin{remark}\quad
\begin{enumerate}\itemsep=0pt
\item Higher order expansions
of these solutions are uniquely determined recursively by the Hamiltonian system and do not
contain any other parameter than the complex parameters $\{\alpha_0, \alpha_1, \alpha_2, \alpha_{\infty},\nu,\eta \}$.

\item These solutions are convergent by Briot--Bouquet's theorem (see Appendix~\ref{appendix2}).

\item The values of complex parameters are generic and should be excluded
the values with which the denominator of the
coef\/f\/icients become zero in the solutions above, that is,
$\{\alpha_1,\alpha_{\infty},\alpha_1\pm\alpha_{\infty},\alpha_0\pm\alpha_2,\alpha_0\pm\alpha_2\pm\alpha_{\infty}\}
\notin \mathbb{Z}$.
\end{enumerate}
\end{remark}

\section{The linear monodromy}

In this section, we calculate the linear monodromy for the solutions (1) and (5).

\subsection{For the solution (1)}\label{subsection4.1}

\subsubsection{The f\/irst limit}

After substituting the solution~(1) into the linear equation~\eqref{lg2}, we take the limit $(t_1, t_2)\rightarrow(0,0)$.
Hereafter
we call this as the f\/irst limit.
Then the linear equation~\eqref{lg2} becomes
\begin{gather}
\frac{d^2 \psi_1}{dx^2}+\left(\frac{2-\alpha_0-\alpha_2}{x}+\frac{1-\alpha_1}{x-1}- \frac{1}{x-b_0}\right) \frac{d
\psi_1}{dx}
\nonumber
\\
\hphantom{\frac{d^2 \psi_1}{dx^2}}
{}+\biggl[\frac{\nu(\nu+\alpha_{\infty})}{x(x-1)}-k_2\left(\frac{1}{x(x-1)}-\frac{1}{x^2}\right)
+\frac{m_2}{x(x-1)(x-b_0)}\biggr]\psi_1=0,
\label{lm:1}
\end{gather}
where $b_0=\frac{\alpha_{\infty}+\alpha_1}{\alpha_{\infty}}$, $k_2=\nu(1-\alpha_0-\alpha_2-\nu)$ and
$m_2=\frac{-\nu\alpha_1}{\alpha_{\infty}}$.
This is a~Heun's type equation with the Riemann scheme
\begin{gather*}
P\left(
\begin{matrix}
x=0\cdot t_2 & x=1 &x=b_0 & x= \infty &
\\
\begin{matrix}
-\nu
\\
-1+\alpha_0+\alpha_2+\nu
\end{matrix}
&
\begin{matrix}
0
\\
\alpha_1
\end{matrix}
&
\begin{matrix}
0
\\
2
\end{matrix}
&
\begin{matrix}
\nu
\\
\nu+\alpha_{\infty}
\end{matrix}
&
\begin{matrix}
;x
\\
\end{matrix}
\end{matrix}
\right).
\end{gather*}
The general solution of~\eqref{lm:1} is
\begin{gather*}
\psi_1=c_1x^{-\nu}+c_2x^{-1+\alpha_0+\alpha_2+\nu}(x-1)^{\alpha_1},
\qquad
 c_1,c_2\in \mathbb{C}.
\end{gather*}
By taking the f\/irst limit, two regular singular points become conf\/luent as
a~regular singular point~\cite{AAB,JM}.
The linear monodromy around $x=0\cdot t_2$ is obtained as a~conf\/luent one $\widetilde M_{t_2}\widetilde M_0$.

The linear monodromy $\{\widetilde M_{t_2}\widetilde M_0, \widetilde M_1,\widetilde M_{\infty}\}$ of~\eqref{lm:1} is
\begin{gather*}
\widetilde M_{t_2}\widetilde M_0=\left(
\begin{matrix} e^{-2\pi i\nu} & 0
\\
0 & e^{2\pi i (\alpha_0+\alpha_2+\nu)}
\end{matrix}
\right),
\qquad
\widetilde M_1=\left(
\begin{matrix} 1 & 0
\\
0 & e^{2\pi i \alpha_1}
\end{matrix}
\right),
\\
\widetilde M_{\infty}=\left(
\begin{matrix} e^{2\pi i\nu} & 0
\\
0 &e^ {2\pi i (\nu+\alpha_{\infty})}
\end{matrix}
\right),
\qquad
\widetilde M_{\infty}\widetilde M_1\widetilde M_{t_2}\widetilde M_0=I_2.
\end{gather*}

We should
separate the conf\/luent linear monodromy $\widetilde M_{t_2}\widetilde M_0$.

\subsubsection{The second limit}

In this section, we separate the conf\/luent linear monodromy $\widetilde M_{t_2}\widetilde M_0$.
After transforming the linear equation~\eqref{lg2} with $x=t_2 \xi$ and substituting the solution~(1) into~\eqref{lg2},
we take the limit $(t_1, t_2) \rightarrow(0,0)$.
Hereafter we call this as the second limit.
By taking the second limit, $x=0$ and $x=t_2$ are separated and $x=1$ and $x=\infty$ become conf\/luent (see Remark~\ref{remark4}).
Then $\psi_2(\xi)=\psi(t_2\xi,t_1,t_2)$ satisf\/ies the following Gauss hypergeometric equation after taking the limit
$(t_1,t_2)\longrightarrow(0,0)$,
\begin{gather}
\frac{d^2 \psi_2}{d\xi^2}+\left(\frac{1-\alpha_0}{\xi}+\frac{1-\alpha_2}{\xi-1}\right) \frac{d
\psi_2}{d\xi}+\frac{k_2}{\xi(\xi-1)}\psi_2=0
\label{lm:2}
\end{gather}
with the Riemann scheme
\begin{gather*}
P\left(
\begin{matrix}
(x=0) & (x=t_2) & (x=1\cdot \infty) &
\\
\begin{matrix}
\xi=0
\\
0
\\
\alpha_0
\end{matrix}
&
\begin{matrix}
\xi=1
\\
0
\\
\alpha_2
\end{matrix}
&
\begin{matrix}
\xi=\infty
\\
\nu
\\
1-\alpha_0-\alpha_2-\nu
\end{matrix}
&
\begin{matrix}
;\xi
\end{matrix}
\end{matrix}
\right).
\end{gather*}
The linear monodromy $\{M_0, M_{t_2}, M_{\infty}M_1 \}$ of~\eqref{lm:2} is
\begin{gather*}
M_0=\left(
\begin{matrix} 1 & 0
\\
0 & e^{2\pi i \alpha_0}
\end{matrix}
\right),
\qquad
M_{t_2}=C_{01}^{-1}\left(
\begin{matrix} 1 & 0
\\
0 & e^{2\pi i \alpha_2}
\end{matrix}
\right)C_{01},
\\
M_{\infty}M_1=C_{0\infty}^{-1}\left(
\begin{matrix} e^{2\pi i\nu} & 0
\\
0 &e^ {2\pi i (\alpha_1+\nu+\alpha_{\infty})}
\end{matrix}
\right)C_{0\infty},
\qquad
M_{\infty} M_1 M_{t_2} M_0=I_2,
\end{gather*}
where $C_{01}$ and $C_{0\infty}$ are the connection matrices of the
Gauss hypergeometric function (see Appendix~\ref{appendix1}, Lemma~\ref{le:1}).
The linear monodromy $\{M_0, M_{t_2}, M_{\infty}M_1 \}$ of~\eqref{lm:2}
is equivalent to $\{\widetilde M_{t_2}\widetilde M_0, \widetilde M_1,\widetilde M_{\infty}\}$ of~\eqref{lm:1}.
We have
\begin{gather*}
M_1=P^{-1}\widetilde M_1P,
\qquad
M_{\infty}=P^{-1}\widetilde M_{\infty}P,
\qquad
M_{\infty}M_1 =P^{-1}\widetilde M_{\infty}\widetilde M_1P
\end{gather*}
for a~matrix $P\in {\rm GL}(2,\mathbb{C})$.
Therefore, $PC_{0\infty}^{-1}\in {\rm GL}(2,\mathbb{C}) $ is a diagonal matrix.
We have
\begin{gather*}
M_1=C_{0\infty}^{-1}\left(
\begin{matrix} 1 & 0
\\
0 & e^{2\pi i \alpha_1}
\end{matrix}
\right)C_{0\infty},
\qquad
M_{\infty}=C_{0\infty}^{-1}\left(
\begin{matrix} e^{2\pi i\nu} & 0
\\
0 &e^ {2\pi i (\nu+\alpha_{\infty})}
\end{matrix}
\right)C_{0\infty}.
\end{gather*}

\begin{remark}
The formal solution of the linear equation~\eqref{lg2} around $x=1$ has the form
\begin{gather*}
\Psi^{(1)}=\left(I_2+\sum\limits_{k=1}^{\infty}\widehat\Psi_k^{(1)}(x-1)^k\right)(x-1)^{T_1}e^{\frac{T}{x-1}},
%\label{xf:1}
\end{gather*}
where $T=\mathrm{diag}(0, \eta t_1)$, $T_1=\mathrm{diag}(0, \alpha_1)$,
$\widehat\Psi_k^{(1)}=\widehat\Psi_k^{(1)}(\lambda_1, \lambda_2,\mu_1,\mu_2,t_1,t_2,\alpha_j,\eta)$.

The series $\widehat\Psi^{(1)} =I_2+\sum\limits_{k=1}^{\infty}\widehat\Psi_k^{(1)}(x-1)^k$ may be a~divergent series
since $x=1$ is an irregular singularity of the Poincar\'e rank one.

The Stokes regions $\mathcal{\widetilde S}_j$ around $x=1$ are given~by
\begin{gather*}
\mathcal{\widetilde S}_j=\left\{x\in \mathbb{C} \, |\, -\varepsilon+(j-1)\pi<{\rm arg}(x-1)<j\pi+\varepsilon,\, |x-1|<r
\right\},
\end{gather*}
where~$\varepsilon$ and~$r$ are suf\/f\/iciently small.
There exist holomorphic functions $\widetilde \Psi_j(x)$ of~\eqref{lg2} on $\mathcal{\widetilde S}_j $ such that
\begin{gather*}
\widetilde \Psi_j(x)\sim\widehat\Psi^{(1)}
\qquad \text{for}\quad
x\rightarrow 1
\qquad \text{and}
\qquad
\Psi_j(x)=\widetilde \Psi_j(x)(x-1)^{T_1}e^{\frac{T}{x-1}}
\end{gather*}
is a~solution of~\eqref{lg2} on $\mathcal{\widetilde S}_j$.
The Stokes matrix $S_j$ is def\/ined~by
\begin{gather*}
\Psi_{j+1}=\Psi_jS_j.
\end{gather*}
We notice that $\Psi_3=\Psi_1(xe^{-2\pi i})e^{2\pi i T_1}$.

First by taking a~limit $t_1\rightarrow 0$ after substituting the solution (1) into the linear equation~\eqref{lg2},
$x=1$ of~\eqref{lg2} becomes a~regular singular point.
Since the coef\/f\/icients $\widehat\Psi_k^{(1)}$ are f\/inite in the limit $t_1\rightarrow 0$, the formal solution
$\Psi^{(1)}$ exists even in the limit $t_1\rightarrow 0$.
Therefore,
\begin{gather*}
\Psi^{(1)}|_{t_1=0}=\left(I_2+\sum\limits_{k=1}^{\infty}\big[\widehat\Psi_k^{(1)}\big]_{t_1=0}(x-1)^k\right)(x-1)^{T_1}
\end{gather*}
has a~regular singularity at $x=1$.
Thus the Stokes matrix $S_j$ become $I_2$ for $j=1,2$.
Then taking a~limit $t_2\rightarrow 0$ makes two regular singular points $x=1$ and $x=\infty$ conf\/luent~\cite{AAB,JM}.
\end{remark}

\subsubsection{The third limit}

In this section, we have Stokes matrices around the irregular singular point $x=1$ by the transformation of the linear
equation~\eqref{lg2}, which keeps the irregularity at $x=1$.
Put
\begin{gather*}
x-1=\frac{\eta t_1}{z},
\end{gather*}
then $\psi_3(z)=\psi(\frac{\eta t_1}{z}+1,t_1,t_2)$ satisf\/ies the following degenerate Kummer's equation after taking
the limit $(t_1,t_2)\rightarrow (0,0)$,
\begin{gather*}
\frac{d^2 \psi_3}{dz^2}+\left(\frac{1+\alpha_1}{z}-\frac{1}{z-\alpha_1}-1\right) \frac{d \psi_3}{dz}=0.
%\label{lm:3}
\end{gather*}
We have the general solution
\begin{gather*}
\psi_3=c_3e^zz^{-\alpha_1}+c_4, \qquad c_3,c_4\in \mathbb{C}.
\end{gather*}
This means the formal solution around the irregular singular point $x=1$ $(z=\infty)$ becomes convergent and Stokes
matrices around $x=1$ become $I_2$.

Therefore, we have the linear monodromy for the solution (1) explicitly.

\begin{theorem}
For the solution $(1)$, the linear monodromy of~\eqref{lg2} is
\begin{gather*}
M_0=\left(
\begin{matrix} 1 & 0
\\
0 & e^{2\pi i \alpha_0}
\end{matrix}
\right),
\qquad
M_{t_2}=C_{01}^{-1}\left(
\begin{matrix} 1 & 0
\\
0 & e^{2\pi i \alpha_2}
\end{matrix}
\right)C_{01},
\\
M_1=S_1^{(1)}S_2^{(1)}C_{0\infty}^{-1}\left(
\begin{matrix} 1 & 0
\\
0 & e^{2\pi i \alpha_1}
\end{matrix}
\right)C_{0\infty},
\qquad
S_1^{(1)}=S_2^{(1)}=I_2,
\\
M_{\infty}=C_{0\infty}^{-1}\left(
\begin{matrix} e^{2\pi i\nu} & 0
\\
0 &e^ {2\pi i (\nu+\alpha_{\infty})}
\end{matrix}
\right)C_{0\infty},
\\
C_{01}=\left(
\begin{matrix} \frac{\Gamma(1-\alpha_0)\Gamma(\alpha_2)} {\Gamma(1-\alpha_0-\nu)\Gamma(\alpha_2+\nu)} &
\frac{\Gamma(1+\alpha_0)\Gamma(\alpha_2)}{\Gamma(1-\nu)\Gamma(\alpha_0+\alpha_2+\nu)}
\vspace{1mm}\\
\frac{\Gamma(1-\alpha_0)\Gamma(-\alpha_2)}{\Gamma(\nu)\Gamma(1-\alpha_0-\alpha_2-\nu)} &
\frac{\Gamma(1+\alpha_0)\Gamma(-\alpha_2)}{\Gamma(\alpha_0+\nu)\Gamma(1-\alpha_2-\nu)}
\end{matrix}
\right),
\\
C_{0\infty}=\left(
\begin{matrix} \frac{e^{\pi i \nu}\Gamma(1-\alpha_0)\Gamma(1-\alpha_0-\alpha_2-2\nu)}
{\Gamma(1-\alpha_0-\alpha_2-\nu)\Gamma(1-\alpha_0-\nu)} & \frac{e^{\pi
i(\alpha_0+\nu)}\Gamma(1+\alpha_0)\Gamma(1-\alpha_0-\alpha_2-2\nu)} {\Gamma(1-\nu)\Gamma(1-\alpha_2-\nu)}
\vspace{1mm}\\
\frac{e^{\pi i(1-\alpha_0-\alpha_2-\nu)}
\Gamma(1-\alpha_0)\Gamma(\alpha_0+\alpha_2+2\nu-1)}{\Gamma(\nu)\Gamma(\alpha_2+\nu)} & \frac{e^{\pi
i(1-\alpha_2-\nu)}\Gamma(1+\alpha_0)\Gamma(\alpha_0+\alpha_2+2\nu-1)}
{\Gamma(\alpha_0+\alpha_2+\nu)\Gamma(\alpha_0+\nu)}
\end{matrix}
\right),
\\
M_{\infty} M_1 M_{t_2} M_0=I_2,
\end{gather*}
for which
$[M_1, M_{\infty}]=0$
holds.
\end{theorem}

\subsection{For the solutions (2), (3) and (4)}

Also for the solutions (2), (3) and (4), we have the similar results, which are summarized as follows.

$\bullet$ For the solution (2):  {\samepage
\begin{gather*}
M_0=\left(
\begin{matrix} 1 & 0
\\
0 & e^{2\pi i \alpha_0}
\end{matrix}
\right),
\qquad
M_{t_2}=C_{01}^{(2)-1}\left(
\begin{matrix} 1 & 0
\\
0 & e^{2\pi i \alpha_2}
\end{matrix}
\right)C_{01}^{(2)},
\\
M_1=S_1^{(1)}S_2^{(1)}C_{0\infty}^{(2)-1}\left(
\begin{matrix} e^{2\pi i \alpha_1} & 0
\\
0 & 1
\end{matrix}
\right)C_{0\infty}^{(2)},
\qquad
S_1^{(1)}=S_2^{(1)}=I_2,
\\
M_{\infty}=C_{0\infty}^{(2)-1}\left(
\begin{matrix} e^{2\pi i\nu} & 0
\\
0 &e^ {2\pi i (\nu+\alpha_{\infty})}
\end{matrix}
\right)C_{0\infty}^{(2)},
\qquad
C_{01}^{(2)}=C_{01}(\nu+\alpha_1, \nu+\alpha_{\infty}, 1-\alpha_0),
\\
C_{0\infty}^{(2)}=C_{0\infty}(\nu+\alpha_1, \nu+\alpha_{\infty}, 1-\alpha_0),
\qquad
M_{\infty} M_1 M_{t_2} M_0=I_2,
\end{gather*}
for which
$
[M_1, M_{\infty}]=0$
holds.}

$\bullet$ For the solution (3):
\begin{gather*}
M_0=\left(
\begin{matrix} 1 & 0
\\
0 & e^{2\pi i \alpha_0}
\end{matrix}
\right),
\qquad
M_{t_2}=\left(
\begin{matrix} 1 & 0
\\
0 & e^{2\pi i \alpha_2}
\end{matrix}
\right),
\\
M_1=S_1^{(1)}S_2^{(1)}C_{01}^{(3)-1}\left(
\begin{matrix} 1 & 0
\\
0 & e^{2\pi i \alpha_1}
\end{matrix}
\right)C_{01}^{(3)},
\qquad
S_1^{(1)}=S_2^{(1)}=I_2,
\\
M_{\infty}=C_{0\infty}^{(3)-1}\left(
\begin{matrix} e^{2\pi i\nu} & 0
\\
0 &e^ {2\pi i (\nu+\alpha_{\infty})}
\end{matrix}
\right)C_{0\infty}^{(3)},
\qquad
C_{01}^{(3)}=C_{01}(\nu, \nu+\alpha_{\infty}, 1-\alpha_0-\alpha_2),
\\
C_{0\infty}^{(3)}=C_{0\infty}(\nu, \nu+\alpha_{\infty},1-\alpha_0-\alpha_2),
\qquad
M_{\infty} M_1 M_{t_2} M_0=I_2,
\end{gather*}
for which
$[M_0, M_{t_2}]=0$
holds.

$\bullet$ For the solution (4):
\begin{gather*}
M_0=\left(
\begin{matrix} 1 & 0
\\
0 & e^{2\pi i \alpha_0}
\end{matrix}
\right),
\qquad
M_{t_2}=\left(
\begin{matrix} 1 & 0
\\
0 & e^{2\pi i \alpha_2}
\end{matrix}
\right),
\\
M_1=S_1^{(1)}S_2^{(1)}C_{01}^{(4)-1}\left(
\begin{matrix} 1 & 0
\\
0 & e^{2\pi i \alpha_1}
\end{matrix}
\right)C_{01}^{(4)},
\qquad
S_1^{(1)}=S_2^{(1)}=I_2,
\\
M_{\infty}=C_{0\infty}^{(4)-1}\left(
\begin{matrix} e^{2\pi i\nu} & 0
\\
0 &e^ {2\pi i (\nu+\alpha_{\infty})}
\end{matrix}
\right)C_{0\infty}^{(4)},
\qquad
C_{01}^{(4)}=C_{01}(\nu,\nu+\alpha_{\infty}, 1-\alpha_2),
\\
C_{0\infty}^{(4)}=C_{0\infty}(\nu, \nu+\alpha_{\infty}, 1-\alpha_2),
\qquad
M_{\infty} M_1 M_{t_2} M_0=I_2,
\end{gather*}
for which
$[M_0, M_{t_2}]=0$
holds.

\subsection{For the solution (5)}

In this section, we calculate the linear monodromy for the solution (5) by the similar way to Subsection~\ref{subsection4.1}.

\subsubsection{The f\/irst limit}

After substituting the solution (5) into the linear equation~\eqref{lg2}, we take the limit $(t_1, t_2)\rightarrow(0,0)$.
At f\/irst we take a~limit $t_1 \rightarrow 0$ keeping $t_2$ as
a~non-zero constant, then we take $t_2\rightarrow 0$.
Then the linear equation~\eqref{lg2} becomes
\begin{gather}
\frac{d^2 \psi_1}{dx^2}+\left(\frac{1-\alpha_0-\alpha_2}{x}+\frac{2-\alpha_1}{x-1}-\frac{1}{x-b_0}\right) \frac{d
\psi_1}{dx}
\nonumber
\\
\hphantom{\frac{d^2 \psi_1}{dx^2}}
{}+\biggl[\frac{\nu(\nu+\alpha_{\infty})}{x(x-1)}-\frac{k_1}{x(x-1)^2} +\frac{m_1}{x(x-1)(x-b_0)}\biggr]\psi_1=0,
\label{lm:4}
\end{gather}
where $b_0=\frac{\alpha_{0}+\alpha_2}{-\alpha_{\infty}}$, $k_1=\nu(\nu+\alpha_1-1)$ and
$m_1=\frac{\nu(\alpha_0+\alpha_2)}{\alpha_{\infty}}$.
This is a~Heun's type equation with the Riemann scheme
\begin{gather*}
P\left(
\begin{matrix}
x=0\cdot t_2 & x=1 &x=b_0 & x= \infty &
\\
\begin{matrix}
0
\\
\alpha_0+\alpha_2
\end{matrix}
&
\begin{matrix}
-\nu
\\
\nu+ \alpha_1-1
\end{matrix}
&
\begin{matrix}
0
\\
2
\end{matrix}
&
\begin{matrix}
\nu
\\
\nu+\alpha_{\infty}
\end{matrix}
&
\begin{matrix}
;x
\\
\end{matrix}
\end{matrix}
\right).
\end{gather*}

The general solution of~\eqref{lm:4} is
\begin{gather*}
\psi_1=c_5(x-1)^{-\nu}+c_6x^{\alpha_0+\alpha_2}(x-1)^{\nu+\alpha_1-1},\qquad c_5,c_6\in \mathbb{C}.
\end{gather*}

The linear monodromy $\{\widetilde M_{t_2}\widetilde M_0, \widetilde M_1,\widetilde M_{\infty}\}$ of~\eqref{lm:4} is
\begin{gather*}
\widetilde M_{t_2}\widetilde M_0=\left(
\begin{matrix} 1 & 0
\\
0 & e^{2\pi i (\alpha_0+\alpha_2)}
\end{matrix}
\right),
\qquad
\widetilde M_1=\left(
\begin{matrix} e^{-2\pi i \nu} & 0
\\
0 & e^{2\pi i (\nu+\alpha_1)}
\end{matrix}
\right),
\\
\widetilde M_{\infty}=\left(
\begin{matrix} e^{2\pi i\nu} & 0
\\
0 &e^ {2\pi i (\nu+\alpha_{\infty})}
\end{matrix}
\right),
\qquad
\widetilde M_{\infty}\widetilde M_1\widetilde M_{t_2}\widetilde M_0=I_2.
\end{gather*}
We should
separate the conf\/luent linear monodromy $\widetilde M_{t_2}\widetilde M_0$.

\subsubsection{The second limit}

In this section, we separate the conf\/luent linear monodromy $\widetilde M_{t_2}\widetilde M_0$.
After transforming the linear equation~\eqref{lg2} with $x=t_2 \xi$ and substituting the solution (5) into~\eqref{lg2},
we take the limit $(t_1, t_2) \rightarrow(0,0)$ as the same as the f\/irst limit.

By taking the second limit, $x=0$ and $x=t_2$ are separated and $x=1$ and $x=\infty$ become conf\/luent.
Then $\psi_2(\xi)=\psi(t_2\xi,t_1,t_2)$ satisf\/ies the following degenerate Heun's equation after taking the limit $(t_1,
t_2) \rightarrow(0,0)$,
\begin{gather}
\frac{d^2 \psi_2}{d\xi^2}+\left(\frac{1-\alpha_0}{\xi}+\frac{1-\alpha_2}{\xi-1}-
\frac{1}{\xi-\xi_{\lambda_2}}\right) \frac{d \psi_2}{d\xi}=0,
\label{lm:5}
\end{gather}
where $\xi_{\lambda_2}=\frac{\alpha_0}{\alpha_0+\alpha_2}$.
We have the general solution
\begin{gather*}
\psi_2=c_7+c_8\xi^{\alpha_0}(\xi-1)^{\alpha_2}, \qquad c_7,c_8 \in \mathbb{C}.
\end{gather*}
The linear monodromy $\{M_0, M_{t_2}, M_{\infty}M_1 \}$ of~\eqref{lm:5} is
\begin{gather*}
M_0=\left(
\begin{matrix} 1 & 0
\\
0 & e^{2\pi i \alpha_0}
\end{matrix}
\right),
\qquad
M_{t_2}=\left(
\begin{matrix} 1 & 0
\\
0 & e^{2\pi i \alpha_2}
\end{matrix}
\right),
\qquad
M_{\infty}M_1=\left(
\begin{matrix} 1 & 0
\\
0 &e^ {-2\pi i (\alpha_0+\alpha_2)}
\end{matrix}
\right),
\\
M_{\infty} M_1 M_{t_2} M_0=I_2.
\end{gather*}

\subsubsection{The third limit}

In this section, we calculate the Stokes matrices around the irregular singular point $x=1$ by the transformation of the
linear equation~\eqref{lg2}, which keeps the irregularity at $x=1$.
Put
\begin{gather*}
x-1=\frac{\eta t_1}{z},
\qquad
\psi=\big(\eta^{-1}z\big)^{\nu}\psi_3(z,t_1,t_2),
\end{gather*}
then $\psi_3(z,t_1,t_2)=(\eta^{-1}z)^{-\nu}\psi(\frac{\eta t_1}{z}+1,t_1,t_2)$ satisf\/ies the following Kummer's
conf\/luent hypergeometric equation after taking the limit $(t_1,t_2)\longrightarrow(0,0)$,
\begin{gather*}
\frac{d^2 \psi_3}{dz^2}+\left(\frac{2\nu+\alpha_1}{z}-1\right)\frac{d \psi_3}{dz}-\frac{\nu}{z}\psi_3=0.
%\label{lm:6}
\end{gather*}
A~system of the fundamental solutions
is
\begin{gather*}
\left({}_1F_1(\nu, 2\nu+\alpha_1; z),
\;
z^{1-2\nu-\alpha_1}{}_1F_1(1-\nu-\alpha_1,2-2\nu-\alpha_1; z) \right).
\end{gather*}
We have the Stokes matrices $S_1^{(1)}$ and $S_2^{(1)}$ around the irregular singular point $x=1$ ($z=\infty$) (see
Appendix~\ref{appendix1}, Lemma~\ref{le:2}).

Summarizing the calculations above, we have the following theorem:

\begin{theorem}%\label{tm:2}
For the solution $(5)$, the linear monodromy of~\eqref{lg2} is
\begin{gather*}
M_0=C^{(5)}\left(
\begin{matrix} 1 & 0
\\
0 & e^{2\pi i \alpha_0}
\end{matrix}
\right)C^{(5)-1},
\qquad
M_{t_2}=C^{(5)}\left(
\begin{matrix} 1 & 0
\\
0 & e^{2\pi i \alpha_2}
\end{matrix}
\right)C^{(5)-1},
\\
M_1=S_1^{(1)}S_2^{(1)} e^{2\pi i T_1},
\qquad
e^{2\pi i T_1}=\left(
\begin{matrix} 1 & 0
\\
0 & e^{2\pi i \alpha_1}
\end{matrix}
\right),
\qquad
S_1^{(1)}=\left(
\begin{matrix} 1 & 0
\\
\frac{-2\pi i e^{\pi i \alpha_1}}{\Gamma(\nu)\Gamma(1-\nu-\alpha_1)} & 1
\end{matrix}
\right),
\\
S_2^{(1)}=\left(
\begin{matrix} 1 &\frac{-2\pi i e^{-2\pi i \alpha_1}}{\Gamma(1-\nu)\Gamma(\nu+\alpha_1)}
\\
0 & 1
\end{matrix}
\right),
\qquad
M_{\infty}=C^{(5)}\left(
\begin{matrix} e^{2\pi i\nu} & 0
\\
0 &e^ {2\pi i (\nu+\alpha_{\infty})}
\end{matrix}
\right)C^{(5)-1},
\\
C^{(5)}=C(\nu, 2\nu+\alpha_1)= \left(
\begin{matrix} \frac{\Gamma(2\nu+\alpha_1)}{\Gamma(\nu+\alpha_1)}e^{\pi i \nu} &
\frac{\Gamma(2-2\nu-\alpha_1)}{\Gamma(1-\nu)}e^{\pi i(1-\nu-\alpha_1)}
\vspace{1mm}\\
\frac{\Gamma(2\nu+\alpha_1)}{\Gamma(\nu)} & \frac{\Gamma(2-2\nu-\alpha_1)}{\Gamma(1-\nu-\alpha_1)}
\end{matrix}
\right),
\\
M_{\infty}M_{t_2} M_0S_1^{(1)}S_2^{(1)}e^{2\pi i T_1}=I_2,
\end{gather*}
for which
$
[M_0, M_{t_2}]=0$,
$[M_0, M_{\infty}]=0$
and
$[M_{t_2}, M_{\infty}]=0$
hold and $C^{(5)}$ is the connection matrix of Kummer's confluent hypergeometric function $($see Appendix~{\rm \ref{appendix1}},
Lemma~{\rm \ref{le:2})}.
\end{theorem}

\begin{remark}\label{remark4}
If we take a~limit $(t_1, t_2) \rightarrow(0,0)$ along a~curve $s_2\sim A t_1^2 (A\in \mathbb{C^{\times}}$), the limit
of the last term in~\eqref{lg2}
\begin{gather*}
\frac{\lambda_2(\lambda_2-1)\mu_2}{x(x-1)(x-\lambda_2)}
\end{gather*}
is not zero.
In our calculation, we take a~special path from $(t_1, t_2)$ to $(0,0)$, such that the numerator of the above term tends
to zero.
Therefore we obtain dif\/ferent limit equations when we choose dif\/ferent paths for the f\/irst and the second limit
equations.
It may be a~contradiction.
But the whole of linear monodromy is the same even though some limit equations are dif\/ferent, since we have the
Riemann--Hilbert correspondence.
In our case, the third limit is the same for any path from $(t_1, t_2)$ to $(0,0)$, which is the main part of the linear monodromy for the solution~(5).
\end{remark}

\subsection{For the  solutions (6), (7) and (8)}
We can determine the linear monodromy for the other solutions (6), (7) and (8), which are summarized as follows.

$\bullet$ For the solution (6):
\begin{gather*}
M_0=C^{(6)}\left(
\begin{matrix} 1 & 0
\\
0 & e^{2\pi i \alpha_0}
\end{matrix}
\right)C^{(6)-1},
\qquad
M_{t_2}=C^{(6)}\left(
\begin{matrix} 1 & 0
\\
0 & e^{2\pi i \alpha_2}
\end{matrix}
\right)C^{(6)-1},
\\
M_{\infty}=C^{(6)}\left(
\begin{matrix} e^{2\pi i(\nu+\alpha_{\infty})} & 0
\\
0 &e^ {2\pi i \nu}
\end{matrix}
\right)C^{(6)-1},
\qquad
M_1=S_1^{(1)}S_2^{(1)}e^{2\pi i T_1},
\\
e^{2\pi i T_1}=\left(
\begin{matrix} 1 & 0
\\
0 & e^{2\pi i \alpha_1}
\end{matrix}
\right),
\qquad
C^{(6)}=C(\nu+\alpha_{\infty}, 2\nu+2\alpha_{\infty}+\alpha_1),
\\
S_1^{(1)}=S_1^{(\infty)}(\nu+\alpha_{\infty}, 2\nu+2\alpha_{\infty}+\alpha_1),
\qquad
S_2^{(1)}=S_2^{(\infty)}(\nu+\alpha_{\infty}, 2\nu+2\alpha_{\infty}+\alpha_1),
\\
M_{\infty}M_{t_2} M_0S_1^{(1)}S_2^{(1)}e^{2\pi i T_1}=I_2,
\end{gather*}
for which
$[M_0, M_{t_2}]=0$,
$[M_0, M_{\infty}]=0$
and
$[M_{t_2}, M_{\infty}]=0$
hold.

$\bullet$ For the solution (7):
\begin{gather*}
M_0=C^{(7)}\left(
\begin{matrix} 1 & 0
\\
0 & e^{2\pi i \alpha_0}
\end{matrix}
\right)C^{(7)-1},
\qquad
M_{t_2}=C^{(7)}\left(
\begin{matrix} e^{2\pi i \alpha_2} & 0
\\
0 & 1
\end{matrix}
\right)C^{(7)-1},
\\
M_{\infty}=C^{(7)}\left(
\begin{matrix} e^{2\pi i\nu} & 0
\\
0 &e^ {2\pi i (\nu+\alpha_{\infty})}
\end{matrix}
\right)C^{(7)-1},
\qquad
M_1=S_1^{(1)}S_2^{(1)}e^{2\pi i T_1},
\\
e^{2\pi i T_1}=\left(
\begin{matrix} 1 & 0
\\
0 & e^{2\pi i \alpha_1}
\end{matrix}
\right),
\qquad
C^{(7)}=C(\nu+\alpha_2, 2\nu+2\alpha_2+\alpha_1),
\\
S_1^{(1)}=S_1^{\infty}(\nu+\alpha_2, 2\nu+2\alpha_2+\alpha_1),
\qquad
S_2^{(1)}=S_2^{\infty}(\nu+\alpha_2, 2\nu+2\alpha_2+\alpha_1),
\\
M_{\infty}M_{t_2} M_0S_1^{(1)}S_2^{(1)}e^{2\pi i T_1}=I_2,
\end{gather*}
for which
$[M_0, M_{t_2}]=0$,
$[M_0, M_{\infty}]=0$
and
$[M_{t_2}, M_{\infty}]=0$
hold.

$\bullet$ For the solution (8):
\begin{gather*}
M_0=C^{(8)}\left(
\begin{matrix} e^{2\pi i \alpha_0} & 0
\\
0 & 1
\end{matrix}
\right)=C^{(8)-1},
\qquad
M_{t_2}=C^{(8)}\left(
\begin{matrix} 1 & 0
\\
0 & e^{2\pi i \alpha_2}
\end{matrix}
\right)C^{(8)-1},
\\
M_{\infty}=C^{(8)}\left(
\begin{matrix} e^{2\pi i\nu} & 0
\\
0 &e^ {2\pi i (\nu+\alpha_{\infty})}
\end{matrix}
\right)C^{(8)-1},
\qquad
M_1=S_1^{(1)}S_2^{(1)}e^{2\pi i T_1},
\\
e^{2\pi i T_1}=\left(
\begin{matrix} 1 & 0
\\
0 & e^{2\pi i \alpha_1}
\end{matrix}
\right),
\qquad
C^{(8)}=C(\nu+\alpha_2+\alpha_{\infty}, 2\nu+2\alpha_2+2\alpha_{\infty}+\alpha_1),
\\
S_1^{(1)}=S_1^{\infty}(\nu+\alpha_2+\alpha_{\infty}, 2\nu+2\alpha_2+2\alpha_{\infty}+\alpha_1),
\\
S_2^{(1)}=S_2^{\infty}(\nu+\alpha_2+\alpha_{\infty}, 2\nu+2\alpha_2+2\alpha_{\infty}+\alpha_1),
\\
M_{\infty}M_{t_2} M_0S_1^{(1)}S_2^{(1)}e^{2\pi i T_1}=I_2,
\end{gather*}
for which
$
[M_0, M_{\infty}]=0$,
$[M_0, M_{t_2}]=0$
and
$[M_{t_2}, M_{\infty}]=0$
hold.

Summarizing the all calculations above, we have the following
theorem:

\begin{theorem}\label{tm:4}
The eight meromorphic solutions around the origin of the two-dimensional degene\-ra\-te Garnier system ${\rm G}_2(1112)$ have the
following characteristics:
\begin{itemize}\itemsep=0pt

\item For the solution $(1)$ and $(2)$, $[M_1, M_{\infty}]=0$ and $S_1^{(1)}=S_2^{(1)}=I_2$ hold.

\item  For the solution $(3)$ and $(4)$, $[M_0, M_{t_2}]=0$ and $S_1^{(1)}=S_2^{(1)}=I_2$ hold.

\item For the solution $(5)$, $(6)$, $(7)$ and $(8)$, $[M_0, M_{\infty}]=0$, $[M_0,M_{t_2}]=0$ and $[M_{t_2},M_{\infty}]=0$ hold.
\end{itemize}
\end{theorem}

\begin{remark}
For the linear monodromy data $\big\{M_0, M_{\infty}, S_1^{(1)}, S_2^{(1)}, e^{2\pi i T_1} \big\}$ of the f\/ifth Painlev\'e
equation, there are three meromorphic solutions around the origin: two solutions such that $[M_0, M_{\infty}]=0$ and one
solution such that $S_1^{(1)}=S_2^{(1)}=I_2$~\cite{KK2}.
\end{remark}

\appendix

\section{Gauss hypergeometric equation and Kummer's equation} \label{appendix1}

In this  appendix,
we show the fundamental solutions and the associated monodromy matrices of Gauss hypergeometric
equation~\cite{KM} and Kummer's equation~\cite{JH}.

\subsection{Gauss hypergeometric equation}

The Gauss hypergeometric equation is
\begin{gather}
x(1-x)\frac{d^2 \psi}{dx^2}+\left(\gamma-(\alpha+\beta+1)x\right)\frac{d\psi}{dx}-\alpha\beta\psi=0.
\label{GH:1}
\end{gather}

 The Riemann scheme of~\eqref{GH:1} is
\begin{gather*}
P \left(
\begin{matrix}
x=0 & x=1 & x= \infty &
\\
\begin{matrix}
0
\\
1-\gamma
\end{matrix}
&
\begin{matrix}
0
\\
\gamma-\alpha-\beta
\end{matrix}
&
\begin{matrix}
\alpha
\\
\beta
\end{matrix}
&
\begin{matrix}
;x
\\
\end{matrix}
\end{matrix}
\right),
\end{gather*}

We list fundamental systems of  solutions for~\eqref{GH:1} around $x=0, 1, \infty$.

$\bullet$ Around $x=0$:
\begin{gather*}
\psi^{(0)} =\left(
\begin{matrix}
{}_2F_1(\alpha, \beta, \gamma; x) & x^{1-\gamma}{}_2F_1(\alpha+1-\gamma, \beta+1-\gamma, 2-\gamma; x)
\end{matrix}
\right).
\end{gather*}

$\bullet$ Around $x=1$:
\begin{gather*}
\psi^{(1)} =  \big(
\begin{matrix}\psi_{1}^{(1)}& \psi_{2}^{(1)}
\end{matrix}
\big),
\\
\psi_{1}^{(1)} = {}_2F_1(\alpha, \beta,\alpha+\beta-\gamma+1;1-x),
\nonumber
\\
\psi_{2}^{(1)} = (1-x)^{\gamma-\alpha-\beta}{}_2F_1(\gamma-\alpha,\gamma-\beta, \gamma+1-\alpha-\beta;1-x).
\nonumber
\end{gather*}

$\bullet$ Around $x=\infty$:
\begin{gather*}
\psi^{(\infty)} =  \big(
\begin{matrix}\psi_{1}^{(\infty)}& \psi_{2}^{(\infty)}
\end{matrix}
\big),
\\
\psi_{1}^{(\infty)} = x^{-\alpha}{}_2F_1\big(\alpha, \alpha-\gamma+1,\alpha+1-\beta;x^{-1}\big),
\nonumber
\\
\psi_{2}^{(\infty)} = x^{-\beta}{}_2F_1\big(\beta, \beta-\gamma+1,\beta+1-\alpha; x^{-1}\big).
\nonumber
\end{gather*}

 The associated monodromy matrices $M_j$ $(j=0, 1, \infty)$ are as follows
\begin{gather*}
M_0=\left(
\begin{matrix} 1 &0
\\
0 & e^{-2 \pi i \gamma}
\end{matrix}
\right),
\qquad
M_1=C_{01}(\alpha, \beta, \gamma)^{-1}\left(
\begin{matrix} 1 &0
\\
0 & e^{2 \pi i (\gamma-\alpha-\beta)}
\end{matrix}
\right)C_{01}(\alpha, \beta, \gamma),
\\
M_{\infty}=C_{0\infty}(\alpha, \beta, \gamma)^{-1}\left(
\begin{matrix} e^{2\pi i \alpha} &0
\\
0 & e^{2 \pi i \beta}
\end{matrix}
\right)C_{0\infty}(\alpha, \beta, \gamma),
\qquad
M_{\infty}M_1M_0=I_2,
\end{gather*}
where $C_{01}(\alpha, \beta, \gamma)$ and $C_{0\infty}(\alpha, \beta, \gamma)$ are connection matrices which are
shown in the following lemma.

\begin{lemma}\label{le:1}
The Gauss hypergeometric function which is the solution of~\eqref{GH:1} has the following connection matrices between
fundamental solutions around two singularities:
\begin{gather*}
\psi^{(i)}=\psi^{(j)}C_{ij}(\alpha, \beta, \gamma), \qquad i,j \in \{0,1,\infty\},
\end{gather*}
where $\psi^{(\nu)}$ $(\nu \in \{0,1,\infty \}) $ is the fundamental solution around the singularity~$\nu$ and
$C_{ij}(\alpha, \beta, \gamma)$ are the connection matrices which are shown as follows
\begin{gather*}
C_{01}(\alpha, \beta, \gamma) = \left(
\begin{matrix} {\Gamma(\gamma)\Gamma(\gamma-\alpha-\beta)}\over {\Gamma(\gamma-\alpha)\Gamma(\gamma-\beta)} & \frac
{\Gamma(2-\gamma)\Gamma(\gamma-\alpha-\beta)} {\Gamma(1-\alpha)\Gamma(1-\beta)}
\vspace{1mm}\\
\frac {\Gamma(\gamma)\Gamma(\alpha+\beta-\gamma)} {\Gamma(\alpha)\Gamma(\beta)}& \frac
{\Gamma(2-\gamma)\Gamma(\alpha+\beta-\gamma)} {\Gamma(1+\alpha-\gamma)\Gamma(1+\beta-\gamma)}
\end{matrix}
\right),
\\
C_{0\infty}(\alpha, \beta, \gamma) = \left(
\begin{matrix} e^{\alpha\pi i} {\Gamma(\gamma)\Gamma(\beta-\alpha)}\over {\Gamma(\beta)\Gamma(\gamma-\alpha)} &
e^{(\alpha-\gamma+1)\pi i} {\Gamma(2-\gamma)\Gamma(\beta-\alpha)}\over {\Gamma(1-\alpha)\Gamma(1-\gamma+\beta)}
\vspace{1mm}\\
e^{\beta \pi i} {\Gamma(\gamma)\Gamma(\alpha-\beta)}\over {\Gamma(\alpha)\Gamma(\gamma-\beta)}& e^{(\beta-\gamma+1)\pi
i} {\Gamma(2-\gamma)\Gamma(\alpha-\beta)}\over {\Gamma(1-\beta)\Gamma(1-\gamma+\alpha)}
\end{matrix}
\right),
\\
C_{\infty1}(\alpha, \beta, \gamma) = \left(
\begin{matrix} {\Gamma(1+\alpha-\beta)\Gamma(\gamma-\alpha-\beta)}\over {\Gamma(\gamma-\beta)\Gamma(1-\beta)} &
{\Gamma(1+\beta-\alpha)} {\Gamma(\gamma-\alpha-\beta)}\over {\Gamma(\gamma-\alpha)\Gamma(1-\alpha)}
\vspace{1mm}\\
e^{(\gamma-\alpha-\beta)\pi i} {\Gamma(1+\alpha-\beta)\Gamma(\alpha+\beta-\gamma)}\over
{\Gamma(1+\alpha-\gamma)\Gamma(\alpha)}& e^{(\gamma-\alpha-\beta)\pi i}
{\Gamma(1+\beta-\alpha)\Gamma(\alpha+\beta-\gamma)}\over {\Gamma(1+\beta-\gamma)\Gamma(\beta)}
\end{matrix}
\right).
\end{gather*}
\end{lemma}

\subsection{Kummer's   equation}

The Kummer's equation is
\begin{gather}
\frac{d^2 \phi}{dx^2}+\left(\frac{\gamma}{x}-1 \right)\frac{d\phi}{dx}-\frac{\alpha}{x}\phi=0.
\label{k:1}
\end{gather}

 The Riemann
scheme of~\eqref{k:1} is
\begin{gather*}
P \left(
\begin{matrix}
   x=0 &  x= \infty & \\
\begin{matrix}
    0\\ 1-\gamma \end{matrix}
& \overbrace{
\begin{matrix}
    0 \\ 1 \end{matrix}
\begin{matrix}
   \alpha \\\quad \gamma-\alpha \end{matrix}}
&\begin{matrix}\quad ;x \\   \end{matrix}
\end{matrix} \right).
\end{gather*}

Fundamental systems of  solutions for~\eqref{k:1} is given by
\begin{gather*}
\phi^{(0)}=\left(
\begin{matrix} {}_1F_1(\alpha, \gamma; x) & x^{1-\gamma}{}_1F_1(\alpha+1-\gamma, 2-\gamma; x)
\end{matrix}
\right).
\end{gather*}

Asymtotic solutions around $x=\infty$ is given by
\begin{gather*}
\phi_1^{(\infty)}\big(e^{-\pi i}x\big)\sim x^{-\alpha}\sum\limits_{k=0}^{\infty}\frac{(-1)^k(\alpha)_k(\alpha+1-\gamma)_k}
{{k!}(e^{-\pi i}x)^k},
\\
\phi_2^{(\infty)}(x)\sim e^x
x^{\alpha-\gamma}\sum\limits_{k=0}^{\infty}\frac{(\gamma-\alpha)_k(1-\alpha)_k}{{k!}x^k}.
\end{gather*}

The associated monodromy matrices are as follows
\begin{gather*}
M_0=C(\alpha, \gamma)\left(
\begin{matrix} 1 &0
\\
0 & e^{-2 \pi i \gamma}
\end{matrix}
\right)C(\alpha, \gamma)^{-1},
\qquad
M_{\infty}=S_1^{(\infty)}(\alpha, \gamma)S_2^{(\infty)}(\alpha, \gamma)e^{2\pi iT_{\infty}},
\\
e^{2\pi iT_{\infty}}=\left(
\begin{matrix} e^{2\pi i \alpha} &0
\\
0 & e^{2 \pi i (\gamma-\alpha)}
\end{matrix}
\right),
\qquad
M_0M_{\infty}=I_2,
\nonumber
\end{gather*}
where $C(\alpha, \gamma)$, $S_1^{(\infty)}(\alpha, \gamma)$ and $S_2^{(\infty)}(\alpha, \gamma)$ are connection matrix and
Stokes matrices respectively
which are
shown in the following lemma.

\begin{lemma}
\label{le:2}
The Kummer's confluent hypergeometric function which is the solution of~\eqref{k:1} has the following connection matrix
between fundamental solutions around $x=0$ and $x=\infty$ and Stokes matrices around $x=\infty$:
\begin{itemize}\itemsep=0pt

\item Connection matrix:
\begin{gather*}
\left(
\begin{matrix}{}_1F_1(\alpha, \gamma; x) & x^{1-\gamma}{}_1F_1(\alpha+1-\gamma, 2-\gamma; x)
\end{matrix}
\right)
=\left(\phi_1^{(\infty)}(e^{-\pi i}x)
\qquad
\phi_2^{(\infty)}(x) \right)C(\alpha, \gamma),
\end{gather*}
where $\phi_i^{(\infty)}$ $(i \in \{1, 2 \}) $ is the fundamental solutions
around the singularity $x=\infty$ and
$C(\alpha, \gamma)$ is the connection matrix which is shown as follows
\begin{gather*}
C(\alpha, \gamma) = \left(
\begin{matrix} \frac{\Gamma(\gamma)e^{\alpha \pi i}} {\Gamma(\gamma-\alpha)} & \frac {\Gamma(2-\gamma)e^{\pi
i(1+\alpha-\gamma)}}{\Gamma(1-\alpha)}
\\
\frac {\Gamma(\gamma)}{\Gamma(\alpha)} & \frac {\Gamma(2-\gamma)}{\Gamma(1+\alpha-\gamma)}
\end{matrix}
\right).
\end{gather*}

\item Stokes matrices:
\begin{gather*}
S_1^{(\infty)}(\alpha, \gamma)= \left(
\begin{matrix} 1 & 0
\\
\frac {-2\pi ie^{\pi i(\gamma-2\alpha)}}{\Gamma(\alpha)\Gamma(1+\alpha-\gamma)} &1
\end{matrix}
\right),
\qquad
S_2^{(\infty)}(\alpha, \gamma)= \left(
\begin{matrix} 1 & \frac {-2\pi ie^{\pi i(4\alpha-2\gamma)}}{\Gamma(1-\alpha)\Gamma(\gamma-\alpha)}
\\
0 &1
\end{matrix}
\right).
\end{gather*}
\end{itemize}
\end{lemma}

\section{Briot--Bouquet's theorem for a~system
of partial dif\/ferential\\ equations in~two variables}
\label{appendix2}

Briot and Bouquet~\cite{BB} showed that existence of a~holomorphic solution for a~special type of nonlinear ordinary
dif\/ferential equations.
In this section we explain the Briot--Bouquet's type theorem for a~system of partial dif\/ferential
equations in two variables following~\cite{GS}.

\subsection{Briot--Bouquet's theorem}

Briot and Bouquet studied a~nonlinear ordinary dif\/ferential equation
\begin{gather}
x\frac{d z}{d x}=h(z,x),
\qquad
z=(z_1,\dots,z_n)
\label{g:0}
\end{gather}
for $h(0,0)=0$.
They have shown
that if the eigenvalues of the Jacobi matrix $(\frac{\partial h}{\partial z}(0,0))$ are not
positive integers, then~\eqref{g:0} has a~convergent holomorphic solution.

R.~Gerard and Y.~Sibuya~\cite{GS}
studied the Briot--Bouquet's type theorem for a~system of partial
dif\/ferential equations in two variables.
They have shown that a~formal solution  will be
convergent:
\begin{lemma}
\label{gs:1}
Assuming that $h_1$ and $h_2$ are holomorphic functions of $z$, $x_1$ and $x_2$ and $z(0,0)=0$, $h_1(0,0,0)=h_2(0,0,0)=0$.
If the simultaneous equations
\begin{gather*}
x_1\frac{\partial z}{\partial x_1}=h_1(z,x_1,x_2),
\qquad
x_2\frac{\partial z}{\partial x_2}=h_2(z,x_1,x_2)
%\label{g:1}
\end{gather*}
have the formal solutions around $(x_1,x_2)=(0,0)$ expressed in power series of $x_1$ and $x_2$, they are convergent.
\end{lemma}

\subsection{Convergence of the solutions}

Solutions (1) and (2) are convergent by Lemma~\ref{gs:1}.
For the solutions with a~pole, for example, solutions (7) and (8), we put
\begin{gather*}
q_1=\frac{Q_1}{t_1},
\qquad
p_1=t_1P_1,
\qquad
q_2=s_2Q_2,
\qquad
p_2= \frac{P_2}{s_2},
\end{gather*}
where $Q_1$, $P_1$, $Q_2$ and $P_2$ are holomorphic functions of $t_1$ and $s_2$ near $(t_1,s_2)=(0,0)$.
Substituting these into the Hamiltonian system $\mathcal{H}_2$, it becomes all Briot--Bouquet's type dif\/ferential
equations with respect to $Q_1$, $P_1$, $Q_2$ and $P_2$.

\subsection*{Acknowledgements}

The author wishes to thank Professor Y.~Ohyama for his constant guidance
and useful suggestions to complete this work.
The author also gives thanks to the anonymous
referees for their relevant contributions to improve this paper.
This work was supported by JSPS KAKENHI Grant Number 22540237 and the Mitsubishi Foundation.

\pdfbookmark[1]{References}{ref}
\LastPageEnding

\end{document}